\documentclass[12pt]{amsart}
\usepackage[english]{babel}
\usepackage{amsfonts,amssymb,latexsym,amscd}

\theoremstyle{plain}
\newtheorem{theorem}{Theorem}

\newtheorem{proposition}{Proposition}
\newtheorem{corollary}{Corollary}
\theoremstyle{definition}
\newtheorem{definition}{Definition}
\theoremstyle{remark}

\usepackage{hyperref}

\begin{document}

\title{Some questions of equivariant movability}
\author{Pavel S. Gevorgyan}
\address{Moscow Pedagogical State University}
\email{pgev@yandex.ru}

\begin{abstract}
In this article some questions of equivariant movability,
connected with the substitution of the acting group $G$  on closed
subgroup $H$ and with transitions to spaces of $H$-orbits and
$H$-fixed points spaces, are investigated. In a special case, the
characterization of equivariantly movable $G$-spaces is given.
\end{abstract}

\keywords{Equivariant shape theory, equivariant movability.}

\subjclass{55P55; 54C56}

\maketitle

\section{Introduction}

This paper is devoted to equivariant movability of $G$-spaces,
i.e., topological spaces endowed with an action of a given compact
group $G$.

More precisely, in \S \ \ref{sec2} we define the notion of
equivariant movability or $G$-movability and we prove several
theorems, including the following ones. If $X$ is $p$-paracompact
and $H \subseteq G$ is a closed subgroup, then $G$-movability of
$X$ implies its $H$-movability (\S \ \ref{sec2}, Theorem
\ref{th1}). $G$-movability of $X$ also implies movability of the
space $X[H]$ of $H$-fixed points in $X$ (\S \ \ref{sec3}, Theorem
\ref{th2}). In particular, equivariant movability of a $G$-space
$X$ implies ordinary movability of the topological space $X$ (\S \
\ref{sec2}, Corollary \ref{cor1}). We construct a non-trivial
example which shows, that the converse, in general, is not true,
even if we take for $G$ the cyclic group $Z_2$ of order 2 (\S \
\ref{sec4}, Example 1). If $X$ is a metrizable $G$-movable space
and $H$ is a closed normal subgroup of $G$, then the space $X|_H$
of its $H$-orbits is also $G$-movable (\S \ \ref{sec5}, Theorem
\ref{th3}). In the case $H=G$ we obtain that $G$-movability of a
metrizable $G$-space implies ordinary movability of the orbit
space $X|_G$ (\S \ \ref{sec5}, Corollary \ref{cor2}). The last
assertion, in general, is not invertible (\S \ \ref{sec5}, Example
2). However, if $X$ is metrizable, $G$ is a compact Lie group and
the action of $G$ on $X$ is free, then $X$ is $G$-movable if and
only if the orbit space $X|_G$ is movable (\S \ \ref{sec6},
Theorem~\ref{th4}). Examples 2 (\S \ \ref{sec5}) and 3 (\S \
\ref{sec7}) show that in the last theorem the assumption that the
group $G$ is a Lie group and the assumption that the action is
free cannot be omitted.

Some of the above listed results with an outline of proof were
given in \cite{gev1}.

Let us denote the category of all topological spaces and
continuous maps by $Top$, the category of all metrizable spaces
and continuous maps by $M$ and the category of all $p$-paracompact
spaces and continuous maps by $P$. Always in this article it is
assumed that all topological spaces are $p$-paracompact spaces and
the group $G$ is compact.

The author is extremely grateful to the referee for his helpful
remarks and comments.

The reader is referred to the books by K. Borsuk \cite{b72} and by
S. Marde\v{s}i\'{c} and J. Segal \cite{ms82} for general
information about shape theory and to the book by G. Bredon
\cite{br72} for introduction to compact transformation groups.

\section{Basic notions and conventions
concerning equivariant topology}\label{sec1}

Let $G$ be a topological group. A topological space $X$ is called
a $G$-space if there is a continuous map $\theta :G \times X \to
X$ of the direct product $G \times X$ into $X$, $\theta(g, x) =
gx$, such that
$$ 1) \quad g(hx)=(gh)x; \qquad 2) \quad ex=x ,$$
for all $g, h \in G$, $x\in X$; here $e$ is the unity of $G$. Such
a (continuous) map $\theta :G \times X \to X$ is called an
(continuous) action of the group $G$ on the topological space $X$.
An evident example is the so called trivial action of $G$ on $X$:
$gx=x$, for all $g\in G, \ x\in X$. Another example is the action
of the group $G$ on itself, defined by $(g, x) \to gx$ for all
$g\in G, \ x\in G$.

If $X$ and $Y$ are $G$-spaces, then so is $X\times Y$, where $g(x,
y)=(gx, gy)$, \  $g\in G, \ (x, y) \in X\times Y$.

A subset $A$ of a $G$ space $X$ is called invariant provided $g\in
G$, $a\in A$ implies $ga\in A$. It is evident, that an invariant
subset of a $G$ space is itself a $G$ space. If $A$ is an
invariant subset of a $G$ space $X$, then every neighborhood of
$A$ contains an open invariant neighborhood of $A$ (see
\cite{p60}, Proposition 1.1.14).

Let $X$ be any $G$-space and let $H$ be a closed and normal
subgroup of the group $G$. The set $Hx=\{hx; \quad h\in H\}$ is
called the $H$-orbit of the point $x \in X$. Clearly the
$H$-orbits of any two points in $X$ are either equal or disjoint,
in other words $X$ is partitioned by its $H$-orbits. We denote the
set of all $H$-orbits of the $G$-space $X$ by $X|_H$. The set
$X|_H$ endowed with the quotient topology is called the $H$-orbit
space of $X$. There is a continuous action of the group $G$ on the
space $X|_H$ defined by the formula $gHx=Hgx, \quad g \in G, x\in
X$. So, $X|_H$ is a $G$-space. In case $H=G$ the $G$-orbit of the
point $x \in X$ is called the orbit of the point $x$ and the
$G$-orbit space is called the orbit space of the $G$-space $X$.

We denote by $X[H]$ the subspace of fixed points of $H$ on $X$, or
the $H$-fixed point subspace of the $G$-space $X$. Let us recall
that $X[H]=\{x \in X; \quad hx=x,$ for any $h \in H\}$.

The set $G_x = \{g\in G; \quad g(x)=x \}$ is a closed subgroup of
the group $G$, for every $x \in X$. $G_x$ is called the stationary
subgroup (or stabilizer) at the point $x$. The action of the group
$G$ on $X$ (or the $G$-space $X$) is called free if the stationary
subgroup $G_x$ is trivial, for every $x\in X$. It is clear that
$G_{gx}=g G_x g^{-1}$, i.e., the stationary subgroups at any two
points of the same orbit are conjugate. The orbits $Gx$ and $Gy$
of points $x$ and $y$, respectively, are said to have the same
type if the stationary subgroups $G_x$ and $G_y$ are conjugate.

Let $X$, \ $Y$ be $G$-spaces. A (continuous) map $f: X\to Y$ is
called a $G$-map, or an equivariant map, if $f(gx)=gf(x)$ for
every $g\in G, \ x\in X$. Note that the identity map $i:X \to X$
is equivariant and the composition of equivariant maps is
equivariant. Therefore, all $G$-spaces and equivariant maps form a
category. Let us denote the category of all topological $G$-spaces
and equivariant maps by $Top_G$, the category of all metrizable
$G$-spaces and equivariant maps by $M_G$ and the category of all
$p$-paracompact $G$-spaces and equivariant maps by $P_G$.

Let $Z$ be a $G$-space and let $Y\subseteq Z$ be an invariant
subset. A $G$-retraction of $Z$ to $Y$ is a $G$-map $r:Z \to Y$
such that $r|_Y=1_Y$.

Let $K_G$ be class of $G$-spaces. A $G$-space $Y$ is called a
$G$-absolute neighborhood retract for the class $K_G$  or a
$G-ANR(K_G)$ ($G$-absolute retract for the class $K_G$  or a
$G-AR(K_G)$), provided $Y\in K_G$ and whenever $Y$ is a closed
invariant subset of a $G$-space $Z\in K_G$, then there exist an
invariant neighborhood $U$ of $Y$ and a $G$-retraction $r:U \to Y$
(there exists a $G$-retraction $r:Z \to Y$).

A $G$-space $Y$ is called a $G$-absolute neighborhood extensor for
the class $K_G$  or a $G-ANE(K_G)$ ($G$-absolute extensor for the
class $K_G$  or a $G-AE(K_G)$), provided for any $G$-space $X\in
K_G$ and any closed invariant subset $A \subseteq X$, every
equivariant map $f:A \to Y$ admits an equivariant extension
$\tilde{f}:U \to Y$, where $U$ is an invariant neighborhood of $A$
in $X$ ($\tilde{f}:X \to Y$).

\section{Movability and equivariant movability}\label{sec2}

The important shape invariant, called movability, was originally
introduced by K. Borsuk \cite{b69} for metric compacta.
Marde\v{s}i\'{c} and Segal \cite{ms70} generalized the notion of
movability to compacta using the $ANR$-system approach. Kozlowski
and Segal in \cite{ks76} gave a categorical description of this
property which applied to arbitrary topological spaces.

Following  Marde\v{s}i\'{c} and Segal \cite{ms70}, let us define
the notion of equivariant movability or $G$-movability :

\begin{definition}
An inverse $G$-system $\underline{X}=\{ X_\alpha , p_{\alpha
\alpha '}, A \}$ where each $X_\alpha $, $\alpha \in A$, is a
$G$-space and every $p_{\alpha \alpha'} : X_{\alpha'} \to X_\alpha
$, $\alpha \leqslant \alpha '$, is a $G$-homotopy class, is called
equivariantly movable or $G$-movable if for every $\alpha \in A$,
there exists an $\alpha ' \in A$, $\alpha '\geqslant \alpha $ such
that for all $\alpha '' \in A$, $\alpha '' \geqslant \alpha $
there exists a $G$-homotopy class $r^{\alpha ' \alpha ''} :
X_{\alpha '} \to X_{\alpha ''}$ such that
$$p_{\alpha \alpha ''} \circ r^{\alpha ' \alpha ''} = p_{\alpha
\alpha '} .$$
\end{definition}

It is known (see \cite{am87}, Theorem 2) that every $G$-space $X$
admits a $G-ANR$-expansion in the sense of Marde\v{s}i\'{c} (see
\cite{ms82}, I, \S \ 2.1), which is the same as saying that there
is an inverse $G-ANR$-system ($G$-system consisting of $G-ANR$'s)
$\underline{X}=\{ X_\alpha , p_{\alpha \alpha '}, A \}$ associated
with $X$ in the sense of Morita \cite{mor75}.

\begin{definition}
A $G$-space $X$ is called equivariantly movable or $G$-movable if
there is an equivariantly movable inverse $G-ANR$-system
$\underline{X}=\{ X_\alpha , p_{\alpha \alpha '}, A \}$ associated
with $X$.
\end{definition}

Note that the last definition of equivariant movability coincides
with the notion of ordinary movability if $G=\{e\}$ is the trivial
group.

Let $X$ be an equivariantly movable $G$-space. The evident
question arises: does movability of the space $X$ follows from its
equivariant movability? The following, more general theorem gives
a positive answer (Corollary 1) to the above question.

\begin{theorem}\label{th1}
 Let $H$ be a closed subgroup of a group $G$. Every $G$-movable $G$-space
 is $H$-movable.
\end{theorem}

To prove this theorem the next result is important.

\begin{theorem}
 Let {H}  be a closed subgroup of a group {G}. Every $G-AR(P_G)$
 $(G-ANR(P_G))$-space is an $H-AR(P_H)$$(H-ANR(P_H))$-space.
\end{theorem}

\begin{proof}
According to a theorem of de Vries (\cite{vr79}, Theorm 4.4), it
is sufficient to show that if $X$ is a $p$-paracompact $H$-space,
then the twisted product $G\times_HX$  is also $p$-paracompact.
Indeed, since $X$ is $p$-paracompact and $G$ is compact, $G\times
X$  is $p$-paracompact. Therefore, the twisted product
$G\times_HX$ is $p$-paracompact.
\end{proof}

\begin{proof}[Proof of Theorem 1]
Let $X$  be any equivariantly movable $G$-space. With respect to
the theorem of Smirnov (\cite{s85}, Theorem 1.3), there is a
closed and equivariant embedding of the $G$-space $X$ to some
$G-AR(P_G)$-space $Y$. Let us consider all open $G$-invariant
neighborhoods of type $F_\sigma$ of the $G$-space $X$ in $Y$. By a
result of R. Palais (\cite{p60}, Proposition 1.1.14), these
neighborhoods form a cofinal family in the set of all open
neighborhoods of $X$ in $Y$, in particular, in the set of all open
and $H$-invariant neighborhoods of the $H$-space $X$ in the
$H$-space $Y$, which, by Theorem 1 is an $H-AR(P_H)$-space. Hence,
from the $G$-movability of the above mentioned family follows its
$H$-movability, i.e. from the $G$-movability of the $G$-space $X$
follows the $H$-movability of the $H$-space $X$.
\end{proof}

From Theorem 1 we obtain the following corollary if we consider
the trivial subgroup $H=\{e\}$ of the group $G$.

\begin{corollary}\label{cor1}
Every equivariantly movable $G$-space $X$ is movable.
\end{corollary}

The converse, in general, is not true, even if one takes for $G$
the cyclic group $Z_2$ of order 2 (see Example 1).

\section{Movability of the $H$-fixed point space}\label{sec3}

\begin{theorem}\label{th2}
Let $H$ be a closed subgroup of a group $G$. If a $G$-space $X$ is
equivariantly movable, then the $H$-fixed point space $X[H]$ is
movable.
\end{theorem}

The proof requires the use of the following theorem.

\begin{theorem}
Let $H$ be a closed subgroup of a group $G$. Let $X$ be a
$G-AR(P_G) (G-ANR(P_G))$- space. Then the $H$-fixed point space
$X[H]$ is an $AR(P) (ANR(P))$-space.
\end{theorem}

\begin{proof}
Let $X$ be a $G-AR(P_G) (G-ANR(P_G))$-space. By Theorem 2, it is
sufficient to prove the theorem in the case $H=G$. I.e., we must
prove that $X[G]$ is $AR(P)$-space. By a theorem of Smirnov
(\cite{s85}, Theorem 1.3), we can consider $X$ as a closed
$G$-subspace of a $G-AR(P_G)$-space $C(G,V)\times \prod {D_\lambda
}$ where $V$ is a normed vector space and thus an $AE(M)$-space,
$C(G,V)$ is the space of continuous maps from $G$ to $V$ with the
compact-open topology and with the action $(g'f)(g)=f(gg'), \quad
g,g' \in G, f \in C(G,V)$ of the group $G$ and $D_\lambda $ is a
closed ball of a finite-dimensional Euclidean space $E_\lambda $
with the orthogonal action of the group $G$.

First, let us prove that the set $(C(G,V)\times \prod {D_\lambda
)}[G]$ of all fixed points of the $G$-space $C(G,V)\times \prod
{D_\lambda }$  is an $AR(P)$-space. The spaces $C(G,V)$ and
$E_\lambda $ are normed spaces. Since the actions of the group $G$
on $C(G,V)$ and $E_\lambda $ are linear, the sets $C(G,V)[G]$ and
$E_\lambda [G]$  will be closed convex sets of locally convex
spaces $C(G,V)$  and $E_\lambda $, respectively. Therefore, by a
well-known theorem of Kuratowski and Dugundji \cite{b67}, $C(G,V)$
and $E_\lambda $ are absolute retracts for metrizable spaces. By a
theorem of Lisica~\cite{l73}, they are also absolute retracts for
$p$-paracompact spaces. For a closed ball $D_\lambda \subset
E_\lambda $ the last conclusion is true  since the set $D_\lambda
[G]=D_\lambda \bigcap E_\lambda [G]$  is closed and convex in
$E_\lambda $.

Since the group $G$ acts on the product $C(G,V)\times \prod
{D_\lambda }$ coordinate-wise, $$ (C(G,V)\times \prod {D_\lambda
)}[G]=C(G,V)[G]\times \left(\prod D_\lambda \right) [G]. $$

Hence, $(C(G,V)\times \prod {D_\lambda )}[G]$ is an $AR(P)$-space,
because it is a product of two $AR(P)$-spaces.

Now let us prove that $X[G]$ is an $AR(P)$-space. Since $X$ is a
$G-AR(P_G)$-space, it is a $G$-retract of the product
$C(G,V)\times \prod {D_\lambda }$. Therefore, $X[G]$ is a retract
of the $AR(P)$-space $(C(G,V)\times \prod {D_\lambda )}[G]$,
hence, it is an $AR(P)$-space.

The absolute neighborhood retract case is proved similarly.
\end{proof}

\begin{proof}[Proof of Theorem 3]
Let $X$ be a $G$-movable space. By Theorem 1, it is sufficient to
prove the theorem in the case $H=G$. So, we must prove movability
of the space $X[G]$ of all $G$-fixed points. We consider the
$G$-space $X$  as a closed and $G$-invariant space of some
$G-AR(P_G)$-space $Y$ (\cite{s85}, Theorem 1.3). The family of all
open, $G$-invariant $F_\sigma$-type neighborhoods $U_\alpha $ of
the $G$-space $X$ in $Y$, is cofinal in the set of all open
neighborhoods of $X$ in $Y$ (\cite{p60}, Proposition 1.1.14). It
consists of $G-ANR(P_G)$-spaces. The intersections $U_\alpha \cap
Y[G] = U_\alpha [G]$ are $ANR(P)$-spaces (Theorem~4). They form a
cofinal family of neighborhoods of the space $X[G]$ in $Y[G]$.
Indeed, for any neighborhood $U$ of the set $X[G]$ in $Y[G]$ there
is a neighborhood $V$ of the set $X[G]$ in $Y$ such that $V \cap
Y[G]=U$. Then the set $W=(Y\setminus Y[G])\cup V$ is a
neighborhood of the set $X$ in $Y$, moreover, $W\cap Y[G]=U$.
There is an $\alpha$ such that $U_\alpha \subset W$ and therefore
$U_\alpha [G] \subset U$. So the family of neighborhoods $U_\alpha
[G]$ is cofinal.

Since $X$ is $G$-movable, for every $U_\alpha $ there is a
neighborhood $U_{\alpha '} \subset U_\alpha$ such that, for any
other neighborhood $U_{\alpha ''} \subset U_{\alpha '}$, there
exists a $G$-equivariant homotopy $F:U_{\alpha '}\times I \to
U_\alpha$ such that $F(y,0)=y$ and $F(y,1)\in U_{\alpha ''}$, for
any $y\in U_{\alpha '}$. It is not difficult to verify that the
homotopy $F[G]:U_{\alpha '}[G]\times I \to U_\alpha [G]$, induced
by $F$,  satisfies the condition of movability of $X[G]$.
\end{proof}

\section{Example of a movable, but not equivariantly movable
space}\label{sec4}

{\bf Example 1.} We will use the idea of S. Marde\v{s}i\'{c}
\cite{m71}. Let us consider the unit circle $S=\{z\in C; \quad
|z|=1\}$. Let us denote $B=[S\times \{1\}]\cup[\{1\}\times S]$.
$B$ is the wedge of two copies of the unit circle $S$ with base
point $\{1\}$. Let us define a continuous map $f:B\to B$ by the
formulas:
\\ $$ f(z,1)=
\begin{cases}
 (z^4, 1),        &\text{$0\leqslant arg(z)\leqslant \frac{\pi}{2}$} \\
 (1, z^4),        &\text{${\frac{\pi}{2}} \leqslant arg(z)\leqslant \pi$} \\
 (z^{-4}, 1),     &\text{$\pi \leqslant arg(z)\leqslant \frac{3\pi}{2}$} \\
 (1, z^{-4}),     &\text{${\frac{3\pi}{2}} \leqslant arg(z) \leqslant 2\pi$}
\end{cases}
$$
$$
f(1,t)=
\begin{cases}
 (t^{-4}, 1),     &\text{$0\leqslant arg(t)\leqslant \frac{\pi}{2}$} \\
 (1, t^{-4}),     &\text{${\frac{\pi}{2}} \leqslant arg(t)\leqslant \pi$} \\
 (t^4, 1),        &\text{$\pi \leqslant arg(t)\leqslant \frac{3\pi}{2}$} \\
 (1, t^4),        &\text{${\frac{3\pi}{2}} \leqslant arg(t) \leqslant 2\pi$}
\end{cases}
$$  \\ for every $z$ and $t$ from $S$. Let us consider the
$ANR$-sequences
\begin{equation*}
B\stackrel{f}{\longleftarrow}
B\stackrel{f}{\longleftarrow}
B\stackrel{f}{\longleftarrow} \cdots
\end{equation*}
and
\begin{equation*}
\Sigma B\stackrel{\Sigma f}{\longleftarrow}
\Sigma B\stackrel{\Sigma f}{\longleftarrow}
\Sigma B\stackrel{\Sigma f}{\longleftarrow} \cdots
\end{equation*}
where $\Sigma $ is the operation of suspension. Let us denote
\begin{equation*}
P=\varprojlim \{B,f\} .
\end{equation*}
Then
\begin{equation*}
\Sigma P=\varprojlim \{\Sigma B,\Sigma f\} .
\end{equation*}
Let us define an action of the group $Z_2=\{e,g\}$ on $\Sigma B$
by the formulas
\begin{equation*}
e[x,t]=[x,t]; \qquad g[x,t]=[x,-t].
\end{equation*}
for every $[x,t] \in \Sigma B, \quad -1 \leqslant t \leqslant 1$. It induces
an action on $\Sigma P$.

\begin{proposition}
The space $\Sigma P$ has trivial shape, but it is not
$Z_2$-movable.
\end{proposition}
\begin{proof}
The triviality of shape of the space $\Sigma P$ is proved by the
method of Marde\v{s}i\'{c} \cite{m71}. Let us prove that the space
$\Sigma P$ is not $Z_2$-movable. Consider the set $\Sigma P[Z_2]$
of all fixed-points of $Z_2$-space $\Sigma P$. It is obvious that
$\Sigma P[Z_2]=P$. Hence, by Theorem 3, it is sufficient to prove
the following proposition.
\end{proof}

\begin{proposition}
The space $P$ is not movable.
\end{proposition}
\begin{proof}
Since the movability of an inverse system remains unchanged under
the action of a functor, it is sufficient to prove non-movability
of the inverse sequence of groups
\begin{equation}
\pi _1(B)\stackrel{f_*}{\longleftarrow} \pi
_1(B)\stackrel{f_*}{\longleftarrow} \pi
_1(B)\stackrel{f_*}{\longleftarrow}\cdots ,
\end{equation}
where $\pi _1(B)$ is the fundamental group of the space $B$ and
$f_*$ is the homomorphism induced by the mapping $f:B\to B$.

It is known that for sequences of groups movability implies the
following condition of Mittag-Leffler, abbreviated as $ML$
(\cite{ms82}, p. 166, Corollary 4):

{\it The inverse system $\{G_\alpha , p_{\alpha \alpha '}, A\}$ of
the $pro-GROUP$ category is said to be $ML$ provided for every
$\alpha \in A$, there exist $\alpha ' \in A, \quad \alpha '
\geqslant \alpha $, such that $p_{\alpha \alpha '}(G_{\alpha
'})=p_{\alpha \alpha ''}(G_{\alpha ''})$, for any $\alpha '' \in
A, \quad \alpha '' \geqslant \alpha $ . }

Thus, it sufficient to prove that the sequence $(1)$ does not
satisfy condition $ML$. Let us observe that $\pi_1(B)$ is a free
group with two generators $a$ and $b$, and $f_*$ is the
homomorphism defined by the formulas
\begin{equation*}
f_*(a)=aba^{-1}b^{-1}, \qquad
f_*(b)=a^{-1}b^{-1}ab.
\end{equation*}

$f_*$ is a monomorphism, because $f_*(a)\neq f_*(b)$, but not an
epimorphism, because, for example, $f_*(x)\neq a$, for all $x\in
\pi_1(B)$. Hence, for any natural $m$ and $n, \quad Imf_*^m
\varsubsetneqq Imf_*^n$ only if $m > n$. It means that the inverse
sequence $(1)$ does not satisfy condition $ML$.
\end{proof}

\section{Movability of the orbit space}\label{sec5}

\begin{theorem}\label{th3}
Let $X$ be a metrizable $G$-space. If $X$ is $G$-movable then for
any closed and normal subgroup $H$ of the group $G$, the $H$-orbit
space $X|_H$ is also $G$-movable.
\end{theorem}

\begin{proof}
Without losing generality one may suppose that $X$ is a closed
$G$-invariant subset of some $G-AR(M_G)$-space $Y$ (\cite{s85},
Theorem 1.1). $X|_H$ is a closed $G$-invariant subset of $Y|_H$
([5], Theorem 3.1).

Let $\{X_\alpha, \alpha \in A \}$ be the family of all
$G$-invariant neighborhoods of $X$ in $Y$. Let us consider the
family $\{X_\alpha |_H, \alpha \in A \}$, where each $X_\alpha |_H
\in G-ANR(M_G)$ and is a $G$-invariant neighborhood of $X|_H$ in
$Y|_H$. Let us prove that the family $\{X_\alpha|_H, \alpha \in A
\}$ is cofinal in the family of all neighborhoods of $X|_H$ in
$Y|_H$. Let $U$ be an arbitrary neighborhood of $X|_H$ in $Y|_H$.
By a theorem of Palais (\cite{p60}, Proposition 1.1.14), there
exists a $G$-invariant neighborhood $V \supset X|_H$ laying in
$U$. Let us denote $\tilde{V} =(pr)^{-1}(V)$, where $pr:Y
\to~Y|_H$ is the $H$-orbit projection. It is evident that
$\tilde{V}$ is a $G$-invariant neighborhood of the space $X$ in
$Y$ and $V=\tilde{V} |_H$. So in any neighborhood of the space
$X|_H$ in $Y|_H$, there is a neighborhood of type $X_\alpha |_H$,
where $X_\alpha $ is a $G$-invariant neighborhood of $X$ in $Y$.

Now let us prove the $G$-movability of the space $X|_H$. Let $X$
be $G$-movable. It means that the inverse system $\{X_\alpha ,
i_{\alpha \alpha '}, A\}$ is $G$-movable. We must prove that the
induced inverse system $\{X_\alpha |_H , i_{\alpha \alpha '}|_H,
A\}$ is $G$-movable. Let $\alpha \in A$ be any index. By the
$G$-movability of the inverse system $\{X_\alpha , i_{\alpha
\alpha '}, A\}$, there is $\alpha ' \in A, \quad \alpha ' >\alpha
$, such that for any other index $\alpha '' \in A, \quad \alpha ''
>\alpha $, there exists a $G$-mapping $r^{\alpha ' \alpha
''}~:~X_{\alpha '}~\to~X_{\alpha ''}$, which makes the following
diagram $G$-homotopy commutative \\

\begin{center}
\begin{picture}(70,70)
\put(0,35){$X_\alpha $}
\put(70,0){$X_{\alpha ''}$}
\put(70,70){$X_{\alpha '}$}

\put(67,71){\vector(-2,-1){54}}
\put(32,63){$i_{\alpha \alpha '}$}
\put(67,5){\vector(-2,1){54}}
\put(33,8){$i_{\alpha \alpha ''}$}
\put(76,66){\vector(0,-1){55}}
\put(78,35){$r^{\alpha ' \alpha ''}$}
\end{picture}
\end{center}
\begin{center}
Diagram 1.
\end{center}

It turns out that, for given $\alpha \in A$, the obtained index
$\alpha ' \in A, \quad \alpha ' > \alpha $, also satisfies the
condition of $G$-movability of the inverse system $\{X_\alpha |_H
, i_{\alpha \alpha '}|_H, A\}$. This is obvious, because the
$G$-homotopy commutativity of Diagram 1 implies the $G$-homotopy
commutativity of the following diagram \\

\begin{center}
\begin{picture}(70,70)
\put(0,35){$X_\alpha|_H $}
\put(70,0){$X_{\alpha ''}|_H$}
\put(70,70){$X_{\alpha '}|_H$}

\put(67,71){\vector(-2,-1){49}}
\put(20,64){$i_{\alpha \alpha '}|_H$}
\put(67,5){\vector(-2,1){49}}
\put(21,7){$i_{\alpha \alpha ''}|_H$}
\put(79,66){\vector(0,-1){55}}
\put(81,35){$r^{\alpha ' \alpha ''}|_H$}
\end{picture}
\end{center}
\begin{center}
Diagram 2. \\
\end{center}
where $r^{\alpha ' \alpha ''}|_H:X_{\alpha '}|_H \to X_{\alpha
''}|_H$ is induced by the mapping $r^{\alpha ' \alpha ''}$. So,
the $G$-movability of the space $X|_H$ is proved.
\end{proof}

\begin{corollary}\label{cor2}
Let $X$ be a metrizable $G$-space. If $X$ is $G$-movable, then the
orbit space $X|_G$ is movable.
\end{corollary}

\begin{proof}
In the case $H=G$ from the last theorem we obtain that the orbit
space $X|_G$ with the trivial action of the group $G$ is
$G$-movable. Therefore, it will be movable by Corollary 1.
\end{proof}

Corollary 2 in general is not invertible: \\ {\bf Example 2.} Let
$\Sigma$ be a solenoid. It is known (\cite{b72}, Theorem 13.5)
that $\Sigma$ is a non-movable compact metrizable Abelian group.
By Corollary 1, the solenoid $\Sigma$ with the natural group
action is not $\Sigma$-movable although the orbit space
$\Sigma|_\Sigma$ as a one-point set is movable .

The converse of Corollary 2 is true if the group $G$ is a Lie
group and the action is free (see Theorem 7).

\section{Equivariant movability of a free $G$-space}\label{sec6}

\begin{theorem}
Let $G$ be a compact Lie group and let $Y$ be a metrizable
$G-AR(M_G)$-space. Suppose that a closed invariant subset $X$ of $
Y$ has an invariant neighborhood whose orbits have the same type.
If the orbit space $X|_G$ is movable, then $X$ is equivariantly
movable.
\end{theorem}

\begin{proof}
The orbit space $X|_G$ is closed in $Y|_G$,  which is a
$G-AR(M)$-space. Let $U$ be an arbitrary invariant neighborhood of
$X$ in $Y$. By the assumption of the theorem, it follows that
there exists a cofinal family of neighborhoods of $X$ in $Y$,
whose orbits have the same type. Therefore, one may suppose that
all orbits of the neighborhood $U$ have the same type. The orbit
set $U|_G$ will be a neighborhood of $X|_G$ in $Y|_G$. From the
movability of $X|_G$ it follows that, for the neighborhood $U|_G$,
there is a neighborhood $\tilde{V}$ of the space $X|_G$ in $Y|_G$,
which lies in the neighborhood $U|_G$ and contracts to any
preassigned neighborhood of the space $X|_G$.

Let us denote $V=(pr)^{-1}(\tilde{V})$, where $pr:Y \to Y|_G$ is
the orbit projection. It is evident that $V$ is an invariant
neighborhood of the space $X$ lying in $U$. Let us prove that $V$
contracts in $U$ to any preassigned invariant neighborhood of $X$.
Let $W$ be any invariant neighborhood of $X$ in $Y$. We must prove
the existence of an equivariant homotopy $F:V \times I \to U$,
which satisfies the condition
\begin{equation*}
F(x,0)=x, \quad F(x,1) \in W ,
\end{equation*}
for any $x \in V$. Since $W|_G$ is a neighborhood of the space
$X|_G$ in $Y|_G$, there is a homotopy $\tilde{F}:V|_G \times I \to
U|_G$ such that
\begin{equation}
F(\tilde{x},0)=\tilde{x}, \quad \tilde{F}(\tilde{x},1) \in W|_G ,
\end{equation}
for any $\tilde{x} \in V|_G$. The homotopy $\tilde{F}:V|_G \times
I \to U|_G$ preserves the $G$-orbit structure, because $V \subset
U$ and all orbits of $U$ have the same types (see Diagram 3).
\\

\begin{center}
\begin{picture}(70,70)
\put(0,0){$V|_G$}
\put(70,0){$U|_G$}
\put(0,70){$V$}
\put(70,70){$U$}

\put(23,3){\vector(1,0){41}} \put(35,5){$i '$}
\put(11,72){\vector(1,0){50}} \put(35,74){$i$}
\put(5,65){\vector(0,-1){53}} \put(-7,35){$pr$}
\put(75,66){\vector(0,-1){55}} \put(77,35){$pr$}
\end{picture}
\end{center}
\begin{center}
Diagram 3. \\
\end{center}
By the covering homotopy theorem of Palais (\cite{p60}, Theorem
2.4.1), there is an equivariant homotopy $F:V\times I \to U$,
which covers the homotopy $\tilde{F}$ and satisfies
$F(x,0)=i(x)=x$. That is, the following diagram is commutative
(Diagram 4). \\

\begin{center}
\begin{picture}(70,70)
\put(0,0){$V|_G \times I$}
\put(70,0){$U|_G$}
\put(0,70){$V \times I$}
\put(70,70){$U$}

\put(42,3){\vector(1,0){22}} \put(48,5){$\tilde{F}$}
\put(32,72){\vector(1,0){30}} \put(43,74){$F$}
\put(10,65){\vector(0,-1){53}} \put(-3,35){$pr$}
\put(75,66){\vector(0,-1){55}} \put(77,35){$pr$}
\end{picture}
\end{center}
\begin{center}
Diagram 4. \\
\end{center}

$F:V \times I \to U$ is the designed  equivariant homotopy. It
only remains to verify that $F(x,1) \in W$. But this immediately
follows from (2) and the commutativity of Diagram 4.
\end{proof}

\begin{theorem}\label{th4}
Let $G$ be a compact Lie group. A metrizable free $G$-space $X$ is
equivariantly movable if and only if the orbit space $X|_G$ is
movable.
\end{theorem}

\begin{proof}
The necessity in a more general case was proved in Corollary 2.
Let us prove the sufficiency. Let the orbit space $X|_G$ be
movable. One can consider the $G$-space $X$ as a closed and
invariant subset of some $G-AR(M_G)$-space $Y$. Let $P \subset X$
be any orbit. From the existence of slices it follows that around
$P$ there is such an invariant neighborhood $U(P)$ in $Y$ that
$typeQ \geqslant typeP$, for any orbit $Q$ from $U(P)$
(\cite{br72}, Corollary 5.5). Since the action of the group $G$ on
$X$ is free, $typeQ=typeP=typeG$, for any orbit $Q$ lying in
$U(P)$. Let us denote $V=\cup \{U(P); \quad P\in X|_G\}$. It is
evident that $V$ is an invariant neighborhood of the space $X$ in
$Y$ and that all of its orbits have the same type. Then, by
Theorem 6, $X$ is equivariantly movable.
\end{proof}

Example 2 shows that the assumption that $G$ is a Lie group is
essential in the above theorem. The Example 3 which follows shows
that the condition of freeness of the action of the group $G$ is
also essential in the above theorem.

\section{Example of a non-free not $Z_2$-movable space with a movable orbit space}\label{sec7}

{\bf Example 3.} Let us consider the space $P=\varprojlim \{B,f\}$
constructed in Example 1. Let us define an action of the group
$Z_2=\{e, g\}$ on the space $B$ by the formulas
\begin{equation}
\aligned e(z,1)&=(z,1) \\ e(1,t)&=(1,t) \\ g(z,1)&=(1,z^{-1}) \\
g(1,t)&=(t^{-1},1) , \\
\endaligned
\end{equation}
for any $z$ and $t$ from $S$. $B$ is a $Z_2-ANR(M_{Z_2})$ space
with the fixed-point $b_0=(1,1)$.

\begin{proposition}
The mapping $f:B \to B$, defined by formulas (3), is equivariant.
\end{proposition}

\begin{proof}
It is necessary to prove the following two equalities:

\begin{equation}
\aligned f(g(z,1))&=g(f(z,1)) \\ f(g(1,t))&=g(f(1,t)) ,\\
\endaligned
\end{equation}
for any $z$ and $t$ from $S$. Let us prove the first one. Consider
the following cases:

{\it Case 1.} $0 \leqslant argz \leqslant \frac{\pi}{2} \quad \Leftrightarrow
\quad \frac{3\pi}{2} \leqslant argz^{-1} \leqslant 2\pi $. \\
Then $f(g(z,1))=f(1,z^{-1})=(1,z^{-4})=g(z^4,1)=gf(z,1)$.

{\it Case 2.} $\frac{\pi}{2} \leqslant argz \leqslant \pi \quad
\Leftrightarrow \quad \pi \leqslant argz^{-1} \leqslant \frac{3\pi}{2}$. \\
Then $f(g(z,1))=f(1,z^{-1})=(z^{-4},1)=g(1,z^4)=gf(z,1)$.

{\it Case 3.} $\pi \leqslant argz \leqslant \frac{3\pi}{2} \quad
\Leftrightarrow \quad \frac{\pi}{2} \leqslant argz^{-1} \leqslant \pi$. \\
Then $f(g(z,1))=f(1,z^{-1})=(1,z^4)=g(z^{-4},1)=gf(z,1)$.

{\it Case 4.} $\frac{3\pi}{2} \leqslant argz \leqslant 2\pi \quad
\Leftrightarrow \quad 0 \leqslant argz^{-1} \leqslant \frac{\pi}{2}$. \\
Then $f(g(z,1))=f(1,z^{-1})=(z^4,1)=g(1,z^{-4})=gf(z,1)$.

The second equality of (4) is proved in a similar way.
\end{proof}

\begin{proposition}
$P$ is a connected, compact, metrizable and equivariantly
non-movable $Z_2$-space which is free at all points except at the
only fixed point $(b_0, b_0, ...)$ and $sh(P|_{Z_2})$=0.
\end{proposition}
\begin{proof}
$P$ is a $Z_2$-space because it is an inverse limit of
$Z_2-ANR(M_{Z_2})$-spaces $B$ and $f$ is an equivariant mapping.
The uniqueness of the fixed point is evident. The connectedness,
compactness and metrizability follows from the properties of
inverse systems (\cite{e77}, Theorem 6.1.20, Corollary 4.2.5). The
non $Z_2$-movability follows from Proposition 2 and Corollary 1.

Let us prove that $sh(P|_{Z_2})=0$ and thus the orbit space
$P|_{Z_2}$ is movable.

Let $X=\varprojlim \{B|_{Z_2},f|_{Z_2}\}$. $X$ is equimorphic to
the orbit space $P|_{Z_2}$. Indeed, let us define a mapping $h:X
\to P|_{Z_2}$ in the following way:
\begin{equation*}
h(([x_1], [x_2], ...))=[(x_1, x_2, ...)]
\end{equation*}
where $([x_1], [x_2], ...) \in X$, and $x_1, x_2, ...$ are
selected from the classes $[x_1], [x_2], ...$ in such way that
$(x_1, x_2, ...) \in P$ or what is the same $f(x_{n+1})=x_n$, for
any $n=1, 2, ... $. Let us prove that the mapping $h$ is defined
correctly. Let $\tilde{x}_1, \tilde{x}_2, ... $ be some other
representatives of the classes $[x_1], [x_2], ... $, respectively,
satisfying the conditions $f(\tilde{x}_{n+1})=\tilde{x}_n$ for any
$n \in N$. Since each class $[x_n]$ has two representatives: $x_n$
and $gx_n$, where $g \in Z_2 = \{e, g\}$, either
$\tilde{x}_n=gx_n$ or $\tilde{x}_n=x_n$. But it is obvious that, if for some $n_0 \in N, \quad \tilde{x}_{n_0}=gx_{n_0}$, then, for any $n \in N, \quad \tilde{x}_n=gx_n$, because $f$ is equivariant. Thus, in the case of another choice of the representatives of the classes $[x_1], [x_2], ... $, we have
\begin{multline*}
h(([x_1],[x_2], ...))=[(\tilde{x}_1,\tilde{x}_2, ...)]=[(gx_1,gx_2, ...)]= \\
=[g(x_1, x_2, ...)]=[(x_1, x_2, ...)] .
\end{multline*}
However, $h$ is a continuous bijection and thus, it is a homeomorphism (\cite{e77}, Theorem 3.1.13).

Consequently,
\begin{equation*}
P|_{Z_2}=\varprojlim\{B|_{Z_2}, f|_{Z_2}\} ,
\end{equation*}
where $B|_{Z_2}\cong S$ and the mapping $\bar{f}=f|_{Z_2}:S \to S$
is defined by the formulas:
\begin{equation}
\bar{f}(z)=
\begin{cases}
 z^4,        &\text{$0\leqslant arg(z)\leqslant \frac{\pi}{2}$} \\
 z^{-4},     &\text{$\frac{\pi}{2}\leqslant arg(z)\leqslant \frac{3\pi}{2}$} \\
 z^4,        &\text{$\frac{3\pi}{2}\leqslant arg(z)\leqslant 2\pi $} \\
\end{cases}
\end{equation}
for any $z \in S$. Thus, we conclude that the orbit space
$P|_{Z_2}$ is a limit of the inverse sequence
\begin{equation*}
S\stackrel{\bar{f}}{\longleftarrow}
S\stackrel{\bar{f}}{\longleftarrow}
S\stackrel{\bar{f}}{\longleftarrow} \cdots
\end{equation*}
By formula (5), the mapping $\bar{f}$ induces a homomorphism
$\bar{f}_*:\pi _1(S) \to \pi _1(S)$, which acts as follows:
\begin{equation*}
\bar{f}_*(a)=aa^{-1}a^{-1}a ,
\end{equation*}
where $a \in \pi _1(S) \cong Z$ is the generator of the group $Z$.
From the above formula, it follows that $\bar{f}_*$ is the null-homomorphism and thus, $deg\bar{f}=0$. For any $k=1, 2, \cdots , \quad \bar{f}_*^k$ is also a null-homomorphism and thus, $deg\bar{f}^k=0$. Therefore, by the classical Hopf theorem (\cite{h59}, Section 2.8, Theorem $H^n$) all $\bar{f}^k:S \to S$ are null-homotopic and $sh(P|_{Z_2})=0$.
\end{proof}

\end{document}